\documentclass[11pt]{amsart}
\usepackage{amsfonts,amsmath,amssymb}
\usepackage{graphicx}
\usepackage{epsfig}

\newcommand{\eps}{\varepsilon}

\newtheorem{lemma}{Lemma}

\newtheorem{theorem}{Theorem}
\newtheorem{corollary}{Corollary}
\newtheorem{example}{Example}
\newtheorem{remark}{Remark}

\newcommand{\R}{\mathbb{R}}
\newcommand{\N}{\mathbb{N}}

\newcommand{\Nc}{\mathcal{N}}
\newcommand{\ONE}{{\bf 1}}

\newcommand{\bpf}[1][Proof]{{\noindent {\sc #1: }}}
\newcommand{\epf}{{{\hspace{4 ex} $\Box$ \smallskip}}}
\renewcommand{\Im}{\mathop{\mathrm{Im}}}

\author{
  Yuri Bakhtin \and Leonid Bunimovich
}
\address{School of Mathematics, Georgia Tech, Atlanta GA, 30332-0160, USA}
\email{bakhtin@math.gatech.edu, bunimovh@math.gatech.edu}
\title{The optimal sink and the best source in a Markov chain}

\begin{document}
\begin{abstract}
It is well known that the distributions of hitting times in Markov chains  are quite irregular, 
unless the limit as time tends to infinity is considered. We show that nevertheless for 
a typical finite irreducible Markov chain  and for  nondegenerate initial distributions the tails 
of the distributions of the hitting times for the states of a Markov chain can be ordered, i.e., 
they do not overlap after a certain finite moment of time.
 If one considers instead each state of a Markov chain as a source rather than a sink then again the 
states can generically be
ordered according to their efficiency. The mechanisms underlying these two orderings are essentially different though.
 Our results can be used, e.g., for a choice of the initial distribution
in numerical experiments with the fastest convergence to equilibrium/stationary distribution, for characterization of the elements of a dynamical network according to their ability to
absorb and transmit the substance (``information'') that is circulated over the network, 
for determining optimal stopping moments (stopping signals/words) when dealing with sequences of symbols, etc.
\end{abstract}
\maketitle

\section{Introduction}
Hitting and recurrence times are a classical subject in the theory of random processes. However,
the relevant studies have always been concerned with averages (expectations) of hitting and recurrence times 
and relations between their distributions for a fixed state \cite{Kac:MR0022323},%
\cite{Galves-Schmitt:MR1483874},\cite{Hirata-Saussol-Vaienti:MR1736991},\cite{Lacroix:MR1952624},%
\cite{Abadi:MR2040782},\cite{Haydn-Vaienti:MR2018869},\cite{Haydn-Lacroix-Vaienti:MR2165587},%
\cite{Kupsa-Lacroix:MR2123204},\cite{Chazottes-Ugalde:MR2151722}.

It is well-known that the distributions of recurrence times are quite regular for many random processes and dynamical systems
\cite{Hirata-Saussol-Vaienti:MR1736991},%
\cite{Haydn-Vaienti:MR2018869},\cite{Haydn-Lacroix-Vaienti:MR2165587},%
\cite{Chazottes-Ugalde:MR2151722}. On the other hand, the distribution functions of the first
hitting times are very irregular 
\cite{Galves-Schmitt:MR1483874},\cite{Lacroix:MR1952624},%
\cite{Abadi:MR2040782},%
\cite{Haydn-Lacroix-Vaienti:MR2165587},%
\cite{Kupsa-Lacroix:MR2123204},\cite{Chazottes-Ugalde:MR2151722},\cite{Keller-Liverani:MR2535206}. This seems to be natural because an ergodic process
returns to any set of positive measure infinitely many times with probability 1, while the hitting event occurs only once.

Therefore, our first result that for a typical irreducible Markov chain and for a typical initial distribution,
the distribution tails of the first hitting times for the states of the chain can be ordered is quite surprising. 
This striking regularity property means that there is
a finite moment of time $n_0$ such that the tails of the survival probabilities $P_i(n)$, $n\ge n_0$, $i=1,2,\ldots,N$, form
an ordered set, i.e., $P_{\sigma_1}(n)<P_{\sigma_2}(n)<\ldots<P_{\sigma_N}(n)$ for all $n\ge n_0$, where
$\sigma_i\in \{1,2,\ldots,N\}$ for all $i$.
From this point of view, the state
$\sigma_1$ is the most efficient sink (absorber of ``information'')  out of all the states of the Markov chain.

The question of the choice of the best (worst) sink naturally arises in the theory of dynamical networks. A dynamical
network is a dynamical system that is generated by individual dynamics of its elements (cells, power stations, neurons, etc.),
the interactions between these elements and the structure of the graph of interactions (often called the topology of the network).
These three characteristics determine the long term dynamics of a network~\cite{Dynamical-networks:MR2335082}.

Traditionally, the theory of dynamical systems deals with asymptotic in time ($t\to\infty$) properties. It has been found
though recently~\cite{Yurchenko} that it is also possible to effectively answer some natural questions on finite time dynamics. 
For instance, placing a hole in a proper place in the phase space of chaotic dynamical systems guarantees that survival
probabilities for this hole for all times $n\ge n_0$ are smaller than for other holes of the same size (measure).

However, it was recently discovered \cite{Yurchenko},\cite{Which-hole:MR2593912} that it is possible  to make finite time predictions of dynamics even if there are no small/large parameters in
equations which govern dynamics of a system.

The results of the present paper (as the ones in \cite{Yurchenko},\cite{Which-hole:MR2593912}) are not only 
generally unexpected but often counterintuitive as well. For instance we provide the examples where the best 
sink or source is not the state with the maximal equilibrium/stationary probability.

It is always tempting and important to try to characterize elements of networks by their ability to absorb and transmit
``information''. By combining the ideas and approaches of \cite{Yurchenko},\cite{Dynamical-networks:MR2335082} it was shown in 
\cite{Which-hole:MR2593912} that, indeed, one can characterize the 
elements of networks by their ability to leak ``information'' out of the system. Thus the elements of networks could
be characterized by their dynamical properties rather than by standard static characteristics like
centrality, betweenness, etc., which are based only on the topology of a network, rather than on its dynamics.

Typically, chaotic dynamical systems, even the most chaotic ones, have a fast decaying but still infinite memory. Therefore, the
studies of statistical properties of dynamical systems always make use of results of the probability theory and, if needed,
require to prove some modifications of the existing limit theorems, etc.  
It is a very natural approach because a memory in chaotic dynamical systems is (most often) infinite and such systems 
are approximated by random processes with a finite memory. However, even for such random processes standard approach 
is to analyze only their asymptotic in time properties. We show here, though, that some interesting finite time properties 
of random processes can also be rigorously studied. For instance, our results show that for hitting times one can find 
not only relations between their averages, but also between their distribution functions. 
It occurred that the infinite tails of these distribution never overlap after a finite moment of time that can be effectively computed.
Our results also generalize those of \cite{Yurchenko},\cite{Which-hole:MR2593912} to 
an essentially larger class of dynamical systems.

However, the question that we address in this paper seems to have never been considered even in the theory of Markov chains.
Our results show that for hitting times one can find not only relations between their averages, but even between their
distributions. 

\smallskip

Another problem considered in this paper is to find the most efficient source in a Markov chain. To the best of our knowledge,
this problem has not been addressed before. It is also motivated by the dynamical networks where the following question is of utmost importance: which node (element of a network) one should apply
a perturbation to, in order to achieve the strongest effect? We show that for a typical irreducible Markov chain, there also exists
a hierarchy of its states with respect to the rate at which the initial perturbation converges to the stationary state.
Thus one can find an optimal node to apply perturbation to in order to achieve the fastest relaxation. 

And again typically there exists a finite moment of time after which the states of a Markov chain
 form an ordered set with respect to their ability to transmit information to the entire chain (network) or to serve as sources.
Generalizations to the case when a sink/source consists not of one but of several states of a Markov chain are straightforward. 
Another straightforward (although important for applications where one deals, e.g., with a network of chemical reactions, supply chains, etc) generalization deals with nonnegative
(rather with transition probabilities) matrices and uses Perron-Frobenius theorem instead of the Markov theorem.
The results of our paper could be used e.g. for choosing an appropriate (e.g. the fastest convergent) initial distribution in computer experiment, for choosing an
appropriate sequence of stopping/observing times when dealing with the sequences of symbols and for dynamical characterization of the elements of networks.

These finite time probabilistic predictions allowed to realize that some natural basic questions have never been addressed in the theory of stochastic processes and even for Markov chains. This gap should be filled in.

\section{Most efficient sink}
It is intuitively clear that for most Markov chains some of the 
states are more important for the dynamics than the others.
The goal of this section is to introduce and study a measure of importance of
the states based on the escape rate through a state (or a family of states, since this generalization of our approach is
straightforward).

Let $P=(P_{ij})_{i,j=1}^N$ be the transition probability matrix of an irreducible Markov chain 
(see, e.g.,~\cite[Chapter XV]{Feller:MR0228020}), on state space
$\{1,\ldots,N\}$ for some $N\in\N$.

Let us fix $k\in\{1,\ldots,N\}$ and stop our Markov chain as soon as it reaches state $k$. In other words,
whenever the original Markov chain makes a transition to $k$, it gets killed, so that the state $k$ can be considered as
a cemetery state for the Markov chain, or a hole through which the mass leaks out of the system.

There are at least two equivalent ways one can describe the resulting dynamics with. One is to treat the new system as a
new Markov
chain with absorbing state $k$ and introduce the associated transition matrix $P^{(k)}$ by
\[
P^{(k)}_{ij}=\begin{cases}
               P_{ij},& i\ne k,\\
               1,& j=i=k,\\
               0,& j\ne i=k.
              \end{cases} 
\]
Another way is to introduce a matrix $Q^{(k)}=(Q_{ij}^{(k)})_{i,j\ne k}$ obtained from $P$ by crossing out its $k$-th row and
column. The matrix $P^{(k)}$ is a stochastic matrix whereas $Q^{(k)}$ is strictly substochastic (or sub-Markov) since
it does not account for the mass leaking out through the state $k$.

We assume that
even after the removal
of an arbitrary state $k$ 
the system remains irreducible and aperiodic, i.e., for some $n_0=n_0(k)$ and all $n>n_0$, all entries of the matrix
$(Q^{(k)})^n$ are positive.
\begin{remark}
 \rm  The aperiodicity assumption is standard, see, e.g.,~\cite[Section XV.9]{Feller:MR0228020}, and we make it
to avoid unnecessary although straightforward technicalities.
\end{remark}

Let us denote the simplex of all probability distributions on $\{1,\ldots,N\}$ by~$\Delta_N$.
Suppose we are given the initial distribution $p=(p_1,\ldots,p_N)\in\Delta_N$. After $n$ steps, the distribution of the Markov
chain with a hole at state $k$ is given by $p(P^{(k)})^n$.
The irreducibility of $P$ implies that, as $n\to\infty$, this distribution converges to the one
concentrated at $k$. This is the only stationary distribution, i.e., the only eigenvector corresponding to the leading
eigenvalue 1 of the stochastic matrix $P^{(k)}$. The rate of convergence to this obvious equilibrium
is characterized by the second largest eigenvalue, $\mu_k<1$. It is easy to see that the spectrum of $P^{(k)}$ coincides
with that of $Q^{(k)}$ except for a simple eigenvalue 1. Therefore, $\mu_k$ is also the leading
positive eigenvalue of matrix $Q^{(k)}$.
In our setting, the
classical Perron--Frobenius 
(PF) theorem guarantees that $\mu_k$ is simple 
and there is an associated eigenvector 
$q^{(k)}=(q^{(k)}_i)_{i\ne k}$ with all positive coordinates.

We can choose $q^{(k)}$ so that
besides the equality
\begin{equation}
\label{eq:q-k-quasistationary}
 q^{(k)} Q^{(k)}=\mu_kq^{(k)},
\end{equation}
it satisfies 
\[\sum_{i\ne k} q^{(k)}_i=1,\]
 thus defining a probability distribution.
Notice that~\eqref{eq:q-k-quasistationary} is exactly the definition of a quasi-stationary distribution for the sub-Markov kernel $Q^{(k)}$.
Since the matrix is sub-Markov, there is no stationary distribution, and the total mass of a vector $qQ^{(k)}$ may
be less than 1 for a probability vector $q$. However, if we normalize the distribution to have total mass 1 after each
step then we end up with the notion of quasi-stationary distributions defined by~\eqref{eq:q-k-quasistationary}. This
equation means that under the stationary distribution, the total mass that has not leaked through $k$
multiplies by $\mu_k<1$ at every step. Therefore, $\lambda_k=-\ln \mu_k$ can serve as the escape rate through $k$. It
can happen that $\mu_k=0$, in this case, all mass escapes the system in finitely many steps, and we set
$\lambda_k=\infty$.

For $q=(q_i)_{i\ne k}$, we define $M_n^{(k)}(q)$ to be the survival probability, or the total mass remaining in the sub-Markov chain defined by $Q^{(k)}$ after
$n$ steps:
\[
 M^{(k)}_{n}(q)=\sum_{i\ne k} (q(Q^{(k)})^n)_i.
\]
If $p=(p_1,\ldots,p_N)$, we define $p^{(k)}$ to be an $N-1$-dimensional vector $(p_i)_{i\ne k}$ and denote
\[
 M^{(k)}_{n}(p)=M^{(k)}_{n}(p^{(k)}).
\]

Since every nonzero vector with nonnegative components has a nontrivial positive component in the direction of the PF eigenvector,
the following statement holds true.

\begin{theorem}\label{th:leak_rate} Let $\lambda_k<\infty$ for some $k\in\{1,\ldots,N\}$.  Then for any $p\in\Delta_N$ with $p_k<1$,
 there are numbers $c_1(p),c_2(p)$ depending only on~$p$ such that
\[
 c_1(p)e^{-\lambda_k n}\le M^{(k)}_{n}(p)\le c_2(p)e^{-\lambda_k n}.
\] 
\end{theorem}

The next corollary compares leaking through different holes.
\begin{corollary} 
\begin{enumerate}
 \item If $\lambda_i>\lambda_j$, then for any $p,q\in\Delta_N$ with $q_j<1$, there is $n_0=n_0(p,q)\in\N$ such that for all $n\ge n_0$,
\[M^{(i)}_{n}(p)< M^{(j)}_{n}(q).\]
\item
Suppose $\sigma$ is a permutation on $\{1,\ldots,N\}$ such that
\[
 \lambda_{\sigma_N}<\ldots<\lambda_{\sigma_1}<\infty.
\]
Then for any family of distributions $(p(i)\in \Delta_N,\ i=1,\ldots,N)$ satisfying $p_i(i)<1$ for all $i$, there is $n_0=n_0(p(1),\ldots,p(N))\in\N$
such that for all $n\ge n_0$,
\[M^{(\sigma_1)}_{n}(p(\sigma_1))< M^{(\sigma_2)}_{n}(p(\sigma_2))<\ldots< M^{(\sigma_N)}_{n}(p(\sigma_N)).\]
\item
 Suppose the state $i\in 1,\ldots,N$ is such that $\lambda_i>\lambda_k$ for all $k\ne i$. For any $p\in\Delta_N$ and any 
$k\ne i$, if $p_k<1$ then there is a time $n_0=n_0(p)$ such that for all $n\ge n_0$,
\[
 M_n^{i}(p)<M_n^k(p).
\]

\end{enumerate}
\end{corollary}
\bpf The first part follows directly from Theorem~\ref{th:leak_rate}. The other two parts are consequences of the first one.\epf

\bigskip

\begin{example}\label{ex:1}\rm Not only the size (stationary probability) of a state of the Markov chain matters for the escape rate through that state.
Consider a Markov chain with transition matrix
 \[ \left( \begin{array}{ccc}
1/3 & 1/6 & 1/2 \\
1/3 & 5/12 & 1/4 \\
1/3 & 5/12 & 1/4 \end{array} \right).
\]
The stationary distribution for this Markov chain is uniform, i.e., $\pi_i=1/3$, $i=1,2,3$. However, the leading eigenvalues in the
reduced matrices $Q^{(i)}$, $i=1,2,3$ are different. Namely $\mu_1=2/3$, $\mu_2=(7+\sqrt{97})/24$, $\mu_3=(9+\sqrt{33})/24$. Therefore,
the fastest escape is through the hole in the third state and the slowest one is through the hole in the second state. This example
belongs to a more general class than the one considered in~\cite{Yurchenko}.
\end{example}

\begin{figure}
\includegraphics[height=6cm]{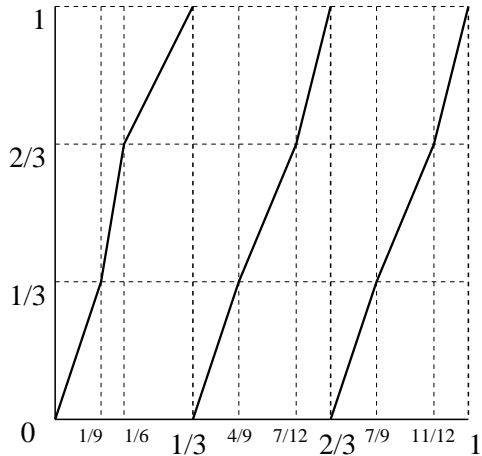}
\caption{Piecewise linear map generating the Markov chain of Example~\ref{ex:1}.}
\label{fig:piecewise}
\end{figure}

\begin{remark}\rm
This Markov chain is generated, e.g., by a piecewise linear map $f:[0,1]\to[0,1]$ shown on Fig.~\ref{fig:piecewise}. States 1,2,3
correspond to intervals $[0,1/3)$, $[1/3,2/3)$, and $[2/3,1]$, respectively, and the stationary measure is Lebesgue measure. 

The Markov chain in the next example is also generated by a certain 1D expanding piecewise linear map. For the sake of brevity we do not present it here though.
\end{remark}

\begin{example}\rm It is possible that the escape is slower through a state with greater stationary probability (bigger ``hole'').
Consider a Markov chain with transition matrix
 \[ \left( \begin{array}{ccc}
1/2 & 1/12 & 5/12 \\
1/2 & 0 & 1/2 \\
1/3 & 1/3 & 1/3 \end{array} \right).
\]
The stationary distribution is given by the vector $(36/83,14/83,33/83)$. The largest eigenvalues of the matrices $Q^{(i)}$,
$i=1,2,3$ equal $\mu_1=(1+\sqrt{7})/6$, $\mu_2=(5+\sqrt{21})/12$, $\mu_3=(3+\sqrt{15})/12$. Therefore, the escape through
the third state is faster than through the first one, although the stationary probability (``size'') of the first state
is larger than that of the third state. This example belongs to a more general class of systems than the one considered in \cite{Which-hole:MR2593912}.
\end{example}

\begin{remark}\rm
 A generalization of Theorem 1 to the case of a non-stochastic but just non-negative matrix is straightforward since the Perron-Frobenius theorem is 
still applicable.
\end{remark}

\begin{remark}\rm  It is easy and straightforward to address the situations where our assumption on irreducibility of the Markov chain
after the removal of a vertex is violated. For example, one can consider the case where besides one strongly connected component $A$
satisfying our original set of assumptions there are 
several extra vertices connected to $A$ but unreachable from $A$. In this case, the rate of escape through any vertex
of $A$ may depend on the initial distribution. In fact, the rate is determined by the minimum 
of the ``internal'' escape rate of $A$ through that vertex and the rates of escape to $A$ from the vertices outside of $A$ that support the 
initial distribution.  Another also easily and directly analyzed situation appears when after a removal of a state the remaining states form
several isolated subsets. Clearly in this case escape from each of these subsets should be 
treated separately (and absolutely analogously to the proof above). It is also easy to see that nothing else besides 
these two situations in case of reducibility of the resulting after removal of the state Markov chain can appear.
\end{remark}

\section{Most efficient source}
In this section we classify the states of a Markov chain with respect to their efficiency in distributing the information or
any perturbation over the entire state space. Here we assume that the Markov chain is just irreducible and aperiodic.

Let the evolution be initiated at state $k\in\{1,\ldots,N\}$. Then for any step $n\ge 0$, the distribution of the
Markov chain at time $n$ is given by~$e_k P^n$
where $e_k$ is the $k$-th coordinate vector and  $P^n$ is the $n$-step transition matrix. We can study the total variation distance between
the distribution at time $n$ and $\pi=(\pi_k)_{k=1}^N$, the stationary distribution, the existence
and uniqueness of which is guaranteed by the Perron--Frobenius Theorem:
\[
 D_k(n)=|e_kP^n-\pi|_1,\quad k\in{1,\ldots,N},\ n\ge 0,
\]
 where $|v|_1=\sum_{i=1}^N|v|$ is the $L^1$ norm of $v$. Ideally, we would like to say that initial state $k_1$ allows for faster convergence to the 
stationary distribution than initial state $k_2$ if there is $n_0\in\N$ such that  $D_{k_1}(n)< D_{k_2}(n)$ for all $n\ge n_0$. However,
there are situations where this property holds due to the specific choice of $|\cdot|_1$ to measure distances, and will be destroyed if
one replaces $|\cdot|_1$ with a different (equivalent) norm. So, we choose to work with a partial order on states that does not depend on the
concrete choice of the norm in $\R^N$.

We denote by $\Nc$ the set of all norms on $\R^N$. We say that a sequence of vectors $(v_n)_{n\in \N}$ in $\R^N$ 
dominates another sequence
of vectors $(u_n)_{n\in \N}$, if for any $H\in\Nc$, there is a number $n_0=n_0(u,v,H)$ such that
\[
H(u_n)<H(v_n),\quad n\ge n_0.
\]
Obviously,  $(v_n)_{n\in \N}$ does not dominate  $(u_n)_{n\in \N}$ iff there is $H\in\Nc$ and a sequence $(n_m)_{m\in\N}$ increasing to infinity such that
\[
 H(u_{n_m})\ge  H(v_{n_m}),\quad m\in\N.
\]

We shall say initial state $k_1$ allows for faster convergence to the 
stationary distribution than initial state $k_2$, if $(e_{k_1}P^n-\pi)_{n\in\N}$ is dominated by $(e_{k_2}P^n-\pi)_{n\in\N}$.

In general, we can say that an initial distribution $u$ allows for faster convergence to the 
stationary distribution than initial distribution $v$, if  $(uP^n-\pi)_{n\in\N}$ is dominated by $(vP^n-\pi)_{n\in\N}$. 
This introduces a partial order on $\Delta_N$, and our goal is to give an equivalent definition of this partial order
in terms of projections on vectors in a (real) canonical Jordan basis $(w_i)_{i=1}^N$ associated to $P$ (we refer to~\cite{HSD-MR2144536}
for the background on canonical forms). 

We assume that $w_N=\pi$, the stationary distribution for $P$, a positive eigenvector of $P$ with simple eigenvalue 1, 
a unique eigenvalue of $P$ equal to~1 in magnitude. 
To each $w_i$, $i=1,\ldots,N-1,$ we associate $\lambda_i$ with $|\lambda_i|<1$ and  $\Im \lambda_i\ge 0$, the eigenvalue of the generalized eigenspace that $w_i$ belongs to
(since complex eigenvalues come in conjugate pairs, we choose $\Im \lambda_i\ge 0$).
 If $\lambda_i\in\R$, then we define
\[r_i=\min\{r\in\N:\ w_i(P-\lambda_i I)^r=0\}.\]
Recalling that for a nonreal eigenvalue $\lambda$, the canonical basis vectors are grouped in pairs, 
we can define $r_i=r_j$ analogously for a pair $(w_i,w_j)$ of canonical basis vectors corresponding to $\lambda_i\not\in\R$. 

In both cases, the numbers $r_i$ enumerate the generalized eigenvectors within one generalized
eigenspace, and the pair $(\lambda_i,r_i)$ determines the rate of decay of $w_i$ under iterations of $P$, namely, $\lambda_i$
is the exponential rate of decay, and $r_i-1$ is the degree of the polynomial factor, see Lemma~\ref{lm:exact-proj} below.

If $\mu\ge 0$ and $k\in\N$, we denote by $\Pi_{\mu,k}v$ the vector projection on the
vector subspace spanned by all $w_i$ such that $|\lambda_i|=\mu$ and $r_i=k$ (the projection
is taken along the span of all other vectors of the Jordan basis). If this subspace is empty, the
projection is assumed to be 0. 

For two vectors $u,v\in \R^N$  we write $u<v$ if there is a real number $a$ with $|a|<1$ such that $u=av$.

\bigskip
\begin{theorem}\label{thm:ordering_initial_distributions} Initial distribution $u$ allows for faster convergence to the 
stationary distribution than initial distribution $v$ if and only if there are $\mu_0\in(0,1)$ and $r_0\in\N$ such that the following
conditions are satisfied:
\begin{enumerate}
 \item[1.] 
If either (i)
$ \mu\in(\mu_0,1)$,
or (ii)
$\mu=\mu_0$ and $r>r_0$, 
then  $\Pi_{\mu,r}u=0$.
\item[2.]
$
\Pi_{\mu_0,r_0} u < \Pi_{\mu_0,r_0} v.$
\end{enumerate} 
\end{theorem}
\begin{remark}\rm Intutively, it is natural to think of $\mu_0$ as of the second largest eigenvalue of $P$. However, the theorem
holds true even in such a degenerate situation where the projections of both $u$ and $v$ on the eigenspace associated to
the second largest eigenvalue vanish. 
\end{remark}

\begin{corollary}\label{cor: second_eigenvalue} 
Suppose that the image of the projection operator $\Pi_{\mu_0,r_0}$ is 1-dimensional  (this is guaranteed if the second largest in magnitude eigenvalue  of $P$ is real and simple).
Let us denote $q_i=|\Pi_{\mu_0,r_0}e_i|$, $i=1,\ldots,N$. Suppose $\sigma$ is a permutation on $\{1,\ldots,N\}$ such that
\[
 q_{\sigma_1}<\ldots<q_{\sigma_N}.
\]
Then for any $N\in\Nc$ there is a number $n_0=n_0(H)$ such that for any $n>n_0$,
\[
 H(e_{\sigma_1}P^n-\pi)<\ldots<H(e_{\sigma_N}P^n-\pi),
\]
so that for any $i,j$ with $i<j$, the initial state $\sigma_i$ allows for faster convergence to the stationary
distribution than the initial state $\sigma_j$. In particular, the initial state $\sigma_1$ allows for faster convergence
than any other initial state.
\end{corollary}

\begin{remark} \rm This hierarchy of states may fail to exist in the case where the dimension of the image of $\Pi_{\mu_0,r_0}$
is greater than one, e.g., where the second highest eigenvalue of $P$ is non-real, or where there are multiple Jordan blocks
associated to $\mu_0$.
\end{remark}

Often, the best source state from the point of view of the hierarchy established in Corollary~\ref{cor: second_eigenvalue}
is the state with maximal stationary probability. However, this is not necessarily so, as the following example shows.

\begin{example}\rm Suppose the transition probability matrix is
 \[ \left( \begin{array}{ccc}
1/8 & 5/8 & 1/4 \\
3/8 & 9/16 & 1/16 \\
1/24 & 1/12 & 7/8 \end{array} \right).
\] 
Then, there are three simple eigenvalues: $1$, $3/4$, and $-3/16$. Their respective eigenvectors are:
$\pi=(1/6,1/3,1/2)$, $w_1=(-1/6,-1/3,1/2)$, and $w_2=(-16/3,13/3,1)$. 
Notice that the stationary probability is maximized by state $3$ since $\pi_3>\pi_2>\pi_1$. However, decomposing
\begin{align*}
e_1&=\pi-\frac{11}{15}w_1-\frac{2}{15}w_2,\\
e_2&=\pi-\frac{17}{15}w_1-\frac{1}{15}w_2,\\
e_3&=\pi+1\cdot w_1+0\cdot w_2,
\end{align*}
comparing the projections on $w_1$, and noticing that $11/15<1<17/15$, 
we can use Theorem~\ref{thm:ordering_initial_distributions} to
conclude that state~1 allows for faster convergence than the two other states. 
\end{example}

\begin{remark}\rm
A generalization of Theorem 2 for non-negative (but non-stochastic) matrices is straightforward. 
\end{remark}

\section{Proof of Theorem~\ref{thm:ordering_initial_distributions}}
We begin with several elementary auxiliary statements. First, we recall formulas for powers of Jordan blocks. For a condition~$A$, we 
use 
\[
\ONE_A=\begin{cases}
        1,& \text{if}\ A\ \text{\rm holds},\\
        0,&\text{\rm otherwise.}
       \end{cases}
\]
\begin{lemma} \label{lm:exact-proj} 1. Let vectors $(w_{i_r})_{r=1}^m$ form a generalized eigenspace of $P$ with eigenvalue
 $\lambda\in\R$, i.e., $w_{i_r}P=\lambda w_{i_r}+\ONE_{2\le r\le m}w_{i_{r-1}}$. Then
 \[
  w_{i_r}P^n=\sum_{k=1}^r 
\binom{n}{r-k}
\lambda^{n-(r-k)}w_{i_k}.
 \]
2. Let $\lambda=\mu e^{i\phi}$, where $\mu>0$ and $\phi\in(0,\pi)$, and vectors $(w_{i_r})_{r=1}^m, (w_{j_r})_{r=1}^m$ form a generalized eigenspace of $P$ with eigenvalue
 $\lambda$, i.e., for any $a,b\in\R$,
\begin{align*}
(aw_{i_r}+bw_{j_r})P
&=\mu (a \cos \phi-b\sin\phi) w_{i_r}+a\ONE_{2\le r\le m}w_{i_{r-1}}\\
 &+ \mu (a \sin \phi+b\cos\phi) w_{j_k}+b\ONE_{2\le r\le m}w_{j_{r-1}}.
\end{align*}
Then
\begin{multline*}
(aw_{i_r}+bw_{j_r})P^n\\
=\sum_{k=1}^r {\binom{n}{r-k}}\mu^{n-(r-k)}\Bigl[(a \cos ((n-(r-k))\phi)-b\sin((n-(r-k))\phi)) w_{i_k}\\
 +(a \sin((n-(r-k))\phi)+b\cos((n-(r-k))\phi)) w_{j_k}\Bigr].
\end{multline*}
\end{lemma}

\begin{lemma}\label{lm:finding_norm} Let $u,v\in\R^N$. If $u<v$, then $H(u)<H(v)$ for any $H\in\Nc$. If $u=v$, then $H(u)=H(v)$ for any $H\in\Nc$. If $u\ne v$ and $u\not<v$, then there is $H\in\Nc$ such that
$H(u)>H(v)$.
\end{lemma}
\bpf First two statements of the lemma are trivial. It is sufficient to prove the third one for the case where $u$ and $v$ are not proportional
to each other. To that end, let us take a linear bijection
that sends vectors $(2,0,0,\ldots,0)$ and $(0,1,0,0,\ldots,0)$ to $u$ and $v$. The pushforward of the Euclidean norm under
this map satisfies the desired property. 
\epf

\begin{lemma}\label{lm:finding_norm_derivative} Suppose $x,y\in\R^N$ and they are not multiples of each other. Then there is $H\in \Nc$,
a neighborhood $U$ of $x$, and a constant $c>0$ such that for all $z\in U$, 
$
 \frac{d}{d\eps} H(z+\eps y)\bigr|_{\eps=0}
$ is well defined and exceeds $c$.
\end{lemma}
\bpf Let us take a linear bijection
that sends vectors $(1,0,0,\ldots,0)$ and $(1,1,0,0,0,\ldots,0)$ to $x$ and $y$. The pushforward of the Euclidean norm under
this map satisfies the desired property. 
\epf

\bpf[Proof of Theorem~\ref{thm:ordering_initial_distributions}]
First, we notice that $\Pi_{1,1}u=\Pi_{1,1}v=\pi=w_N$. This follows from the following facts: $P$-iterates of $u$ and $v$ converge to $\pi$;
$\Pi_{1,1}u$ and $\Pi_{1,1}v$ are invariant under $P$; $P$-iterates of all other projections decay exponentially.

Suppose that there are $\mu_0$ and $r_0$ such that conditions 1 and 2 of the Theorem hold true.
 Decomposing $u$ and $v$ w.r.t.\ the canonical basis and using Lemma~\ref{lm:exact-proj}, we immediately
see that $u$ allows for faster convergence than~$v$.

Suppose now that $u$ allows for faster convergence than~$v$. Let us choose $\mu_0$ and $r_0$ so that  
$\Pi_{\mu_0,r_0}v\ne 0$ and 
if either (i)
$ \mu\in(\mu_0,1)$,
or (ii)
$\mu=\mu_0$ and $r>r_0$, 
then  $\Pi_{\mu,r}v=0$.

Lemma~\ref{lm:exact-proj} immediately implies now that condition 1 of the theorem is satisfied. To prove condition 2, let us
assume that the opposite holds, i.e., $\Pi_{\mu_0,r_0} u\not< \Pi_{\mu_0,r_0} v$. First, we consider the case where $\Pi_{\mu_0,r_0} u \ne \Pi_{\mu_0,r_0} v$.
Lemma~\ref{lm:finding_norm} allows us to find a norm $H\in\Nc$ such that
$H(\Pi_{\mu_0,r_0} u)> H(\Pi_{\mu_0,r_0} v)$. For any small neighborhoods $U$ of $\Pi_{\mu_0,r_0}u$ and $V$ of $\Pi_{\mu_0,r_0}v$ we can use Lemma~\ref{lm:exact-proj} to find an infinite sequence
$n_m\to\infty$ such that
\begin{align}
 \frac{\Pi_{\mu_0,r_0} u P^{n_m}} {\binom{n_m}{r_0 -1}\mu_0^{n_m}}\in U,\quad m\in\N,\label{eq:return_u}
\\ \frac{\Pi_{\mu_0,r_0} v P^{n_m}}{\binom{n_m}{r_0 -1}\mu_0^{n_m}}\in V,\quad m\in\N.\label{eq:return_v}
\end{align}
This is trivially true with $n_m\equiv m$  if all the eigenvalues with magnitude~$\mu_0$ are real and equal to $\mu_0$. If the arguments
of some of these eigenvalues are not zero, then we can use the recurrence property of the shift on the multidimensional
torus induced by these arguments.

Using Lemma~\ref{lm:exact-proj} to compute the leading terms of $(u-\pi) P^{n_m}$  and $(v-\pi) P^{n_m}$, we see that 
\begin{equation}
 \lim_{m\to\infty}\frac{(z-\pi) P^{n_m}-(\Pi_{\mu_0,r_0} z) P^{n_m}}{\binom{n_m}{r_0 -1}\mu_0^{n_m}}=0,\quad z=u,v.
\label{eq:main_contrib_is_mu_0_k_0}
\end{equation}
Therefore, 
\begin{align}
 \frac{(u-\pi) P^{n_m}} {\binom{n_m}{r_0 -1}\mu_0^{n_m}}\in U,\quad m\in\N,\label{eq:return_u_1}
\\ \frac{(v-\pi) P^{n_m}}{\binom{n_m}{r_0 -1}\mu_0^{n_m}}\in V,\quad m\in\N.\label{eq:return_v_1}
\end{align}
so that choosing $U$ and $V$ disjoint and sufficiently small and using inequality $H(\Pi_{\mu_0,r_0} u)> H(\Pi_{\mu_0,r_0} v)$ along with the continuity of $H$,
we conclude that $u$ does not allow for faster convergence than $v$. This contradicts our assumption and therefore it remains to consider
\begin{equation}
\Pi_{\mu_0,r_0} u = \Pi_{\mu_0,r_0} v.
\label{eq:equal_at_mu_0_k_0}
\end{equation}
Faster convergence for $u$ is clearly impossible in the situation where $u=v$. Assuming $u\ne v$, we can find
numbers $\mu_1$ and $r_1$ such that $\Pi_{\mu_1,r_1}u\ne \Pi_{\mu_1,r_1}v$ and if (i) $\mu_1<\mu<\mu_0$, or (ii) $\mu=\mu_1$ and $r>r_1$, then
$\Pi_{\mu,r}u=\Pi_{\mu,r}v$.

Let $U$,$H$, and $c$ be defined in Lemma~\ref{lm:finding_norm_derivative} applied to $x=\Pi_{\mu_0,k_0} u = \Pi_{\mu_0,k_0} v$,
and $y=\Pi_{\mu_1,r_1} u-\Pi_{\mu_1,r_1} v$ which is not a multiple of $x$.
Due to Lemma~\ref{lm:exact-proj}, there is a sequence of numbers $n_m\to\infty$ and a sequence of vectors
$u_m,v_m$
such that
\begin{align}
\label{eq:decomposition1}
 \frac{ (u-\pi) P^{n_m}}{ \binom{n_m}{r_0 -1}\mu_0^{n_m}}=z_m+u_m,\quad m\in\N,\\
\label{eq:decomposition2} \frac{(v-\pi) P^{n_m}}{\binom{n_m}{r_0 -1}\mu_0^{n_m}}=z_m+v_m,\quad m\in\N,
\end{align}
where $z_m\in U$ for all $m$,  and
\begin{align}
u_m=\frac{\binom{n_m}{r_1 -1}\mu_1^{n_m} \Pi_{\mu_1,r_1} u}{\binom{n_m}{r_0 -1}\mu_0^{n_m}}
+o\left(\frac{\binom{n_m}{r_1 -1}\mu_1^{n_m}}{\binom{n_m}{r_0 -1}\mu_0^{n_m}}\right),\quad m\to\infty,
\label{eq:decomposition3}\\
v_m=\frac{\binom{n_m}{r_1 -1}\mu_1^{n_m} \Pi_{\mu_1,r_1} v}{\binom{n_m}{r_0 -1}\mu_0^{n_m}}+o\left(\frac{\binom{n_m}{r_1 -1}\mu_1^{n_m}}{
\binom{n_m}{r_0 -1}\mu_0^{n_m}}\right),\quad m\to\infty,
\label{eq:decomposition4}
\end{align}
We can use relations \eqref{eq:decomposition3} and \eqref{eq:decomposition4} to derive
\begin{multline*}
H(z_m+u_m)-H(z_m+v_m)=\frac{\binom{n_m}{r_1 -1}\mu_1^{n_m}}{\binom{n_m}{r_0 -1}\mu_0^{n_m}}\cdot \frac{d}{d\eps}H(z_m+\eps (\Pi_{\mu_1,r_1} u - \Pi_{\mu_1,r_1} v))\bigr|_{\eps=0}
\\+o\left(\frac{\binom{n_m}{r_1 -1}\mu_1^{n_m}}{\binom{n_m}{r_0 -1}\mu_0^{n_m}}\right),\quad m\to\infty.
\end{multline*}
Since $z_m\in U$, 
 Lemma~\ref{lm:finding_norm_derivative} allows us to conclude that the derivative in the r.h.s.\ of the last identity exceeds $c>0$.
Therefore, relations \eqref{eq:decomposition1} and \eqref{eq:decomposition2} imply that
$H((u-\pi)P^{n_m})>H((v-\pi)P^{n_m})$ for all $m$, which contradicts our assumption that $u$ allows for faster convergence than $v$.
Therefore,~\eqref{eq:equal_at_mu_0_k_0} is impossible and the proof of the necessity of conditions 1 and 2 of the theorem is complete.
\epf

\section{Concluding remarks}
We have shown that the tails of the distributions of hitting times for different states of irreducible Markov chains and for typical initial 
distributions can be ordered. This means that there is a finite moment of time $n^*$ after which the tails of these
distributions never intersect. This property allows to determine the optimal sink in a Markov chain or in a dynamical
network.

Our results hold for any nondegenerate initial distribution and in this respect they essentially generalize 
those in~\cite{Which-hole:MR2593912}, where only Lebesgue measure was considered.

We also demonstrated that one can determine the best source in a Markov chain.  
Again it is a finite time result and the hierarchy of the Markov chain states emerges in their ability to serve as a source. 
For a network, this suggests the node or element one should apply a perturbation to, or inject information at, so that the perturbation spreads over
the network and converges to the stationary distribution in the fastest way.
Our results are also true (with obvious adjustments) if the matrix $P$ is nonnegative and not necessarily stochastic.

\section{Acknowledgments}
Y.B.~was partially supported by NSF CAREER grant DMS-0742424. L.B.~was partially supported by NSF grant DMS-0900945.

\bibliographystyle{plain}
\bibliography{mc}

\end{document}